# Stiff-PINN: Physics-Informed Neural Network for Stiff Chemical Kinetics

Weiqi Ji [a†*], Weilun Qiu [b†], Zhiyu Shi [b†], Shaowu Pan [c], Sili Deng [a*]

[a] *Department of Mechanical Engineering, Massachusetts Institute of Technology,*

*Cambridge, MA 02139, USA*

[b] *College of Engineering, Peking University, Beijing 100871, China*

[c] *Department of Aerospace Engineering, University of Michigan, Ann Arbor, MI 48109, USA*

[†] These authors contributed equally to this work.

[*] Corresponding Author: Sili Deng (silideng@mit.edu)

## Abstract

Recently developed physics-informed neural network (PINN) has achieved success in many science and engineering disciplines by encoding physics laws into the loss functions of the neural network, such that the network not only conforms to the measurements, initial and boundary conditions but also satisfies the governing equations. This work first investigates the performance of PINN in solving stiff chemical kinetic problems with governing equations of stiff ordinary differential equations (ODEs). The results elucidate the challenges of utilizing PINN in stiff ODE systems. Consequently, we employ Quasi-Steady-State-Assumptions (QSSA) to reduce the stiffness of the ODE systems, and the PINN then can be successfully applied to the converted non/mild-stiff systems. Therefore, the results suggest that stiffness could be the major reason for the failure of the regular PINN in the studied stiff chemical kinetic systems. The developed Stiff-PINN approach that utilizes QSSA to enable PINN to solve stiff chemical kinetics shall open the possibility of applying PINN to various reaction-diffusion systems involving stiff dynamics.

State Assumptions.

## 1. Introduction

Deep learning has enabled advances in many scientific and engineering disciplines, such as computer visions, natural language processing, and autonomous driving. Depending on the applications, many different neural network architectures have been developed, including Deep Neural Networks (DNN), Convolutional Neural Networks (CNN), Recurrent Neural Networks (RNN), and Graph Neural Network (GNN). Some of them have also been employed for data-driven physics modeling [1–8], including turbulent flow modeling [9] and chemical kinetic modeling [10–14]. Those different neural network architectures introduce specific regularization to the neural network based on the nature of the task such as the scale and rotation invariant of the convolutional kernel in CNN. Among them, the recently developed Physics-Informed Neural Network approach (PINN) [15–21] enables the construction of the solution space of differential equations using deep neural networks with space and time coordinates as the inputs. The governing equations (mainly differential equations) are enforced by minimizing the residual loss function using automatic differentiation and thus it becomes a physics regularization of the deep neural network. This framework permits solving differential equations (i.e., forward problems) and conducting parameter inference from observations (i.e., inverse problems). PINN has been employed for predicting the solutions for the Burgers' equation, the Navier–Stokes equations, and the Schrodinger equation [16]. To enhance the robustness and generality of PINN, multiple variations of PINN have also been developed, such as Variational PINNs [22], Parareal PINNs [23], and nonlocal PINN [24].

Despite the successful demonstration of PINN in many of the above works, Wang et al. [25] investigated a fundamental mode of failure of PINN that is related to numerical stiffness leading to unbalanced back-propagated gradients between the loss function of initial/boundary conditions and the loss function of residuals of the differential equations during model training. In addition to the numerical stiffness, physical stiffness might also impose new challenges in the training of PINN. While PINN has been applied for solving chemical reaction systems involving a single-step reaction [19], stiffness usually results from the nonlinearity and complexity of the reaction network, where the characteristic time scales for species span a wide range of magnitude. Consequently, the challenges for PINN to accommodate stiff kinetics can potentially arise from several reasons, including the high dimensionality of the state variables (i.e., the number of species), the high nonlinearity resulted from the interactions among species, the imbalance in the loss functions for different state variables since the species concentrations could span several orders of magnitudes. Nonetheless, stiff chemical kinetics is essential for the modeling of almost every real-world chemical system such as atmospheric chemistry and the environment, energy conversion and storage, materials and chemical engineering, biomedical and pharmaceutical engineering. Enabling PINN for handling stiff kinetics will open the possibilities of using PINN to facilitate the design and optimization of these wide ranges of chemical systems.

In chemical kinetics, the evolution of the species concentrations can be described as ordinary differential equation (ODE) systems with the net production rates of the species as the source terms. If the characteristic time scales for species span a wide range of magnitude, integrating the entire ODE systems becomes computationally intensive. Quasi-Steady-State-Assumptions (QSSA) have been widely adopted to simplify and solve stiff kinetic problems, especially in the 1960s when efficient ODE integrators were unavailable [26]. A canonical example of the utilization of QSSA is the Michaelis–Menten kinetic formula, which is still widely adopted to formulate enzyme reactions in biochemistry. Nowadays, QSSA is still widely employed in numerical simulations of reaction-transport systems to remove chemical stiffness and enable the explicit time integration with relatively large time steps [27–29]. Moreover, imposing QSSA also reduces the number of state variables and transport equations by eliminating the fast species such that the computational cost can be greatly reduced. From a physical perspective [26,30], QSSA identifies the species (termed as QSS species) that are usually radicals with relatively low concentrations. Their net production rates are much lower than their consumption and production rates and thus can be assumed zero. From a mathematical perspective [26], the stiffness of the ODEs can be characterized by the largest absolute eigenvalues of the Jacobian matrix, i.e., the Jacobian matrix of the reaction source term to the species concentrations. QSSA identifies the species that correspond to the relatively large eigenvalues of the chemical Jacobian matrix and then approximate the ODEs with differential-algebraic equations to reduce the magnitude of the largest eigenvalue of the Jacobian matrix and thus the stiffness.

In the current work, we will evaluate the performance of PINN in solving two classical stiff dynamics problems and compare it with the performance of Stiff-PINN, which incorporates QSSA into PINN to reduce stiffness. In Section 2, the two classical stiff kinetic systems, the corresponding PINN models, and the implementation of QSSA to formulate the Stiff-PINN models will be presented. In Section 3, the performances of the regular-PINN and Stiff-PINN in solving the two stiff problems will be investigated. Finally, conclusions and the outlook of future work will be presented.

## 2. Methodology

### 2.1 Stiff Chemical Kinetic Systems

A homogenous chemical reaction system can be modeled using the following ordinary differential equations (ODEs):

$$\frac{d\boldsymbol{y}}{dt} = \boldsymbol{f}(t, \boldsymbol{y}), \quad t_0 \leq t \leq t_{final} \tag{2.1}$$

$$\boldsymbol{y}(t_0) = \boldsymbol{y}_0, \tag{2.2}$$

where $\boldsymbol{y} = [y_1, y_2, \ldots, y_N]^T$ is the column vector of species concentrations, and $N$ is the number of chemical species. $t$ is the time, and the initial and final time are denoted as $t_0$ and $t_{final}$, respectively. $\boldsymbol{y}_0$ is the column vector of the initial species concentrations. The ODE system described by Eq. (2.1) with the initial conditions specified by Eq. (2.2) can be numerically solved using an ODE integrator such as the explicit Euler method or Runge-Kutta method. However, many ODEs for chemical kinetic models are stiff [31], and solving stiff ODEs with an explicit method requires very small time steps such that the integration could be computationally intensive. Otherwise, implicit ODE integrators such as Backward Differentiation Formula can be used. However, in general, solving stiff ODEs is time-consuming since the implicit method usually involves solving the nonlinear systems with Newton iteration. Therefore, it is still an active research area to efficiently solve stiff ODE systems [32,33], which is an integral part of many reaction-diffusion systems, such as in chemical engineering, energy conversion, and biomedical applications.

While it is difficult to give a precise definition of the stiffness of a chemical kinetic model, one criterion can be whether there are largely separated time scales for different species. For instance, some of the fast-evolving species have very short time scales while some of the species evolve very slowly and have orders of magnitude longer time scales. To resolve those species concentrations with short time scales $\tau_{fast}$, one has to use very small time steps in the explicit ODE integrators. However, to resolve the slowly evolving species, the number of integration steps scale with $S = \tau_{final}/\tau_{fast}$. If $S$ is on the order of 1000 or larger, the system will be considered as stiff [31]. However, the shortest time scale $\tau_{fast}$ is defined locally and is evolving, such that it is difficult to define the stiffness of a problem precisely. A practical approach to measure the stiffness is to compare the computational cost of explicit ODE integrators developed for non-stiff problems and implicit ones for stiff problems on a specific problem. If the computational cost using an implicit ODE integrator is much lower than the explicit ODE integrators, the problem can be regarded as a very stiff problem.

This work will investigate the performance of PINN in two classical stiff chemical kinetic problems, ROBER [34] and POLLU [35], which are extensively used for testing stiff ODE integrators. Specifically, the ROBER problem [34] consists of three species and five reactions, and the POLLU problem [35] consists of 20 species and 25 reactions describing the air pollution formation in atmospheric chemistry. The formula of the ROBER problem (three ODEs) is presented here to illustrate the implementation of PINN in the following sections, while the formula of the POLLU problem (20 ODEs) is presented in the Supporting Information.

The ROBER problem refers to the following reaction network,

$$A \xrightarrow{k_1} B,$$
$$B + B \xrightarrow{k_2} C + B, \qquad (2.3)$$
$$B + C \xrightarrow{k_3} A + C.$$

The reaction rate constants are $k_1 = 0.04, k_2 = 3 \times 10^7, k_3 = 10^4$, and the initial conditions are $y_1(0) = 1$, $y_2(0) = 0$, $y_3(0) = 0$, where $y_1, y_2, y_3$ denote the concentrations of A, B, C, respectively. The evolution of the species concentrations can be described by the following ODEs:

$$\begin{aligned} \frac{dy_1}{dt} &= -k_1 y_1 + k_3 y_2 y_3, \\ \frac{dy_2}{dt} &= k_1 y_1 - k_2 y_2^2 - k_3 y_2 y_3, \\ \frac{dy_3}{dt} &= k_2 y_2^2. \end{aligned} \qquad (2.4)$$

The reaction rate constants vary in a range of nine orders of magnitude, i.e., $k_2/k_1 \sim 10^9$, resulting in a system with strong stiffness.

## 2.2 Physics-Informed Neural Network

Without loss of generality, we shall use the ROBER problem to illustrate the framework of the PINN. Figure 1 illustrates the structure of PINN informed by the ODEs of the ROBER problem.

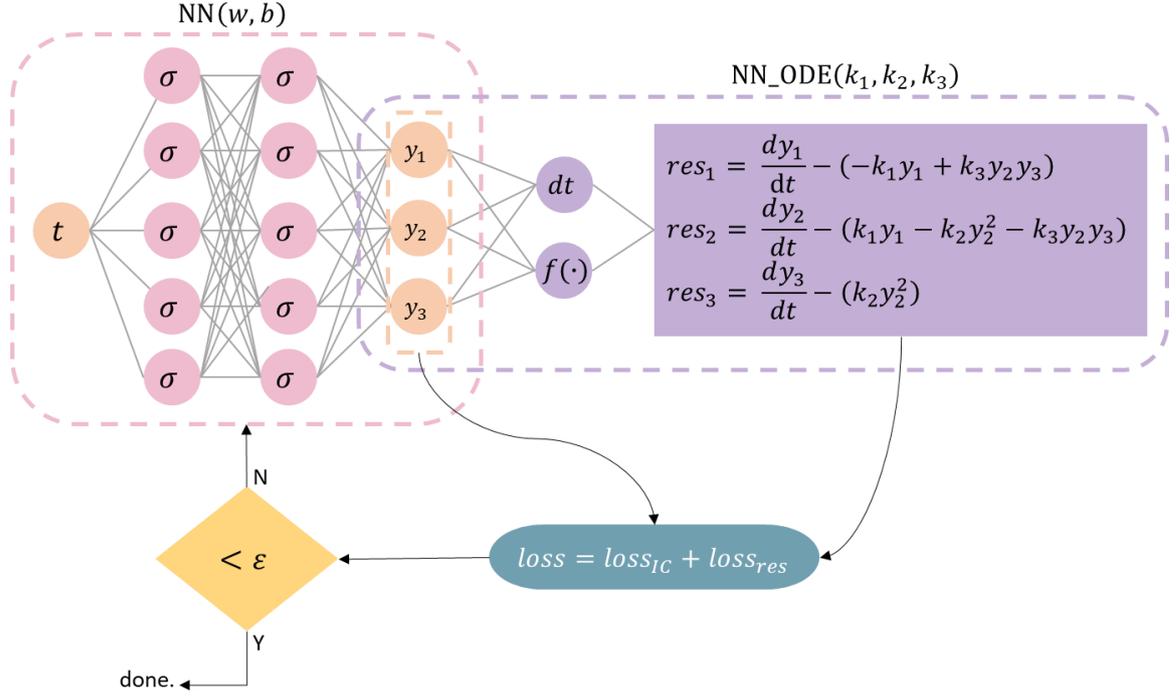

Figure 1. A schematic of the regular-PINN for the ROBER problem. The only input is the time $t$, and the output is the solution vector $[y_1(t), y_2(t), y_3(t)]^T$, which has to satisfy the governing equations and the initial conditions. The two neutral networks, NN($w$, $b$) and NN_ODE($k_1$, $k_2$, $k_3$), share parameters, and both contribute to the loss function.

The PINN framework consists of two neural networks. The first part of the framework is a neural network NN($w$, $b$) that takes the time $t$ as the input and outputs the concentrations of all of the species $\boldsymbol{y} = [y_1, y_2, y_3]^T$ at that time. Then the output of NN is fed into a second network NN_ODE($k_1$, $k_2$, $k_3$), which is essentially the governing differential equations of the ROBER problem, to evaluate the residuals of the ODEs. Finally, a loss function is constructed by combining the loss functions of the initial conditions and residuals. Specifically,

$$loss = loss_{IC} + loss_{res},$$
$$loss_{IC} = \sum_{i=1}^{3} w_{IC}^i [y_i^{NN}(t=t_0) - y_i(t=t_0)],$$
$$loss_{res} = \frac{1}{R}\sum_{j=1}^{R}\sum_{i=1}^{3} w_{res}^i [res_i(t=t_j)]. \qquad (2.5)$$

where the weights of the initial conditions and the residuals for different species in the loss functions are rescaled by $w_{IC}^i$ and $w_{res}^i$, since they could be unbalanced with each other by orders of magnitude and hinder the training. The predicted initial conditions to evaluate $loss_{IC}$ are from the first neural network NN, and the loss of the residuals is from the second network NN_ODE. The residual loss $loss_{res}$ is evaluated at randomly sampled points in the computational domains, i.e., $\{t_1, t_2, \ldots, t_R\} \in [t_0, t_{final}]$. Backpropagation through the two networks is conducted on the auto-differentiation framework of PyTorch to compute the gradient of the loss functions to the weights of the neural networks. The time derivatives of $dy_i/dt$ used in the residual loss is obtained using auto-differentiation as well. The neural network is optimized via the first-order optimizer Adam [36].

**2.3 Quasi-Steady-State Assumptions**

As previously discussed, the QSSA is often imposed on certain species to reduce the computational cost of simulating the evolution of the system. These QSS species, which are often radicals and unstable intermediates, generally have shorter time scales compared to other species. By assuming that the net production rates of the QSS species are zero, the concentrations of these species can be expressed by algebra equations instead of ODEs, such that the number of ODEs to model the kinetic system is reduced. Note that although the net production rates of the QSS species are small, the production and consumption rates of them are not necessarily small. In a general form, the net production rate of species $k$ can be written as

$$\frac{dY_k}{dt} = \omega_k^+ - \omega_k^-. \tag{2.6}$$

$Y_k$ is the species concentration, and $\omega_k^+, \omega_k^-$ are the production and consumption rate, respectively. QSSA requires that

$$\left|\frac{dY_k}{dt}\right| \ll (\omega_k^+, \omega_k^-), \tag{2.7}$$

so that

$$\omega_k^+ - \omega_k^- \approx 0. \tag{2.8}$$

Take the ROBER problem as an example, as shown in Fig. 2, the concentration of species $y_2$ increases sharply during the initial induction period, and then it goes to a phase of slow change with considerably low concentrations compared to $y_1$ and $y_3$. The production rate of $y_2$ is comparable with the consumption rate of $y_1$, as can be seen from Eq. (2.3). However, since the net production rate over the entire integration range can be estimated to be proportional to the maximum species concentrations and the maximum concentration of $y_2$ is five orders of magnitude lower than that of $y_1$, the net production rate of $y_2$ is much slower than $y_1$. Similarly, the consumption rate of $y_2$ is comparable with the production rate of $y_3$, while the net production rate of $y_2$ is much slower than $y_3$. Therefore, the species $y_2$ is likely to be a QSS species compared to $y_1$ and $y_3$.

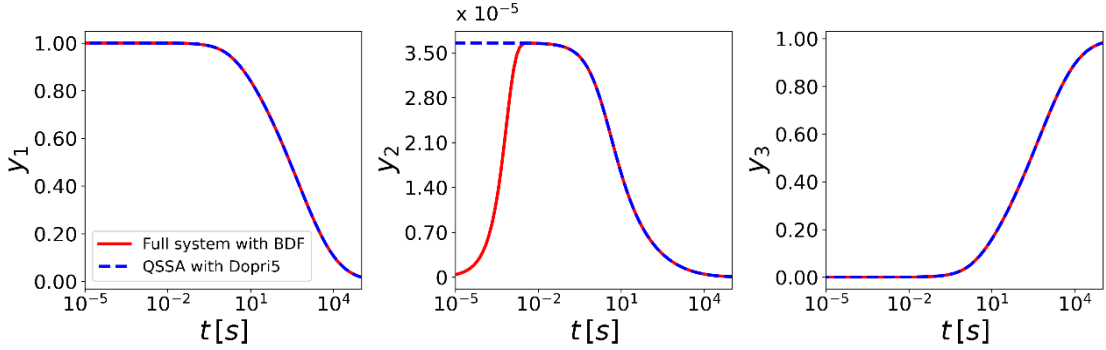

Figure 2. The comparisons of the solutions of the ROBER problem using BDF solver and the solutions of the reduced system via QSSA using Dopri5. QSSA does not affect the solution of $y_1$ and $y_3$ and accurately predicts $y_2$ after the initial induction period. Note that the time is presented in the logarithmic scale to better illustrate the evolution of $y_2$ during the induction period.

The assumption of $y_2$ as a QSS species implies that

$$0 = k_1 y_1 - k_2 y_2^2 - k_3 y_2 y_3 \tag{2.9}$$

and

$$y_2 = \frac{-k_3 y_3 + \sqrt{k_3^2 y_3^2 + 4 k_1 k_2 y_1}}{2 k_2}. \qquad (2.10)$$

Consequently, the original ROBER problem of three ODEs can be approximated by the following differential-algebraic equations (DAEs) and two ODEs.

$$\begin{aligned} \frac{dy_1}{dt} &= -k_1 y_1 + k_3 y_2 y_3, \\ y_2 &= \frac{-k_3 y_3 + \sqrt{k_3^2 y_3^2 + 4 k_1 k_2 y_1}}{2 k_2}, \\ \frac{dy_3}{dt} &= k_2 y_2^2. \end{aligned} \qquad (2.11)$$

The issue of stiffness in the system described by Eq. (2.11) should be reduced compared to the system described by Eq. (2.4) since the fast-evolving species $y_2$ is not explicitly solved. Therefore, we can use an explicit integrator to solve Eq. (2.10) and the results are shown in Fig. 2. It is shown that the explicit Runge-Kutta method of Dopri5 method well predicts $y_1$ and $y_3$ under QSSA. In addition, the species profile of $y_2$ can be computed based on the algebraic equation in Eq. (2.10), and the approximated profile agrees well with the accurate solution of $y_2$ except for the initial induction period. Note that, in practice, what interests us are the stable species rather than unstable intermediates.

We can then integrate QSSA into the framework of PINN, and the new framework denoted as Stiff-PINN is illustrated in Fig. 3. The key differences between stiff-PINN and regular-PINN are three folds: first, NN_ODE is informed by the reduced systems with QSSA, instead of the original ODEs. Second, NN will only output the non-QSS species, and the QSS species are approximated in NN_ODE. Third, the QSS species are excluded from the loss functions of the initial conditions and the loss functions of the residuals, and therefore, loss functions will only consist of the information of non-QSS species.

It is worthy to note that manually deriving explicit algebraic expressions for the QSS species could be challenging for complex systems. Automatic model reduction tools [37,38] could help tackle the challenges by linearizing the QSSA approximation. Alternatively, the QSS species have to be solved via an additional nonlinear optimization step.

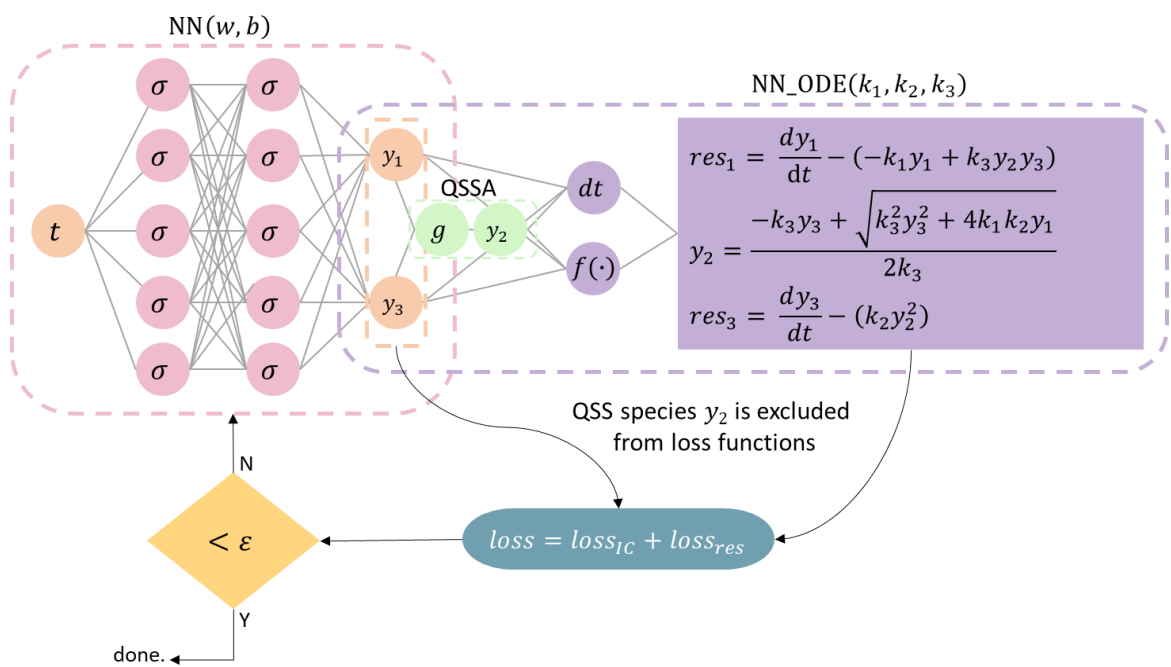

Figure 3. A schematic of the Stiff-PINN for the ROBER problem. Comparing the model pipeline for Stiff-PINN with that for the regular-PINN, species $y_2$ is assumed in quasi-state state (QSS) and the residuals of these QSS species are excluded from the loss functions during the training of Stiff-PINN. By removing these relatively fast-evolving species from the loss function, the stiffness of the differential equations is reduced.

## 3. Results and Discussions

In this section, the performance of regular-PINN in two classical stiff chemical kinetic problems: ROBER and POLLU will be investigated. As will be discussed, the stiffness is the main reason for the failure of regular-PINN, and the performance of Stiff-PINN in removing stiffness and predicting the species profiles will then be demonstrated. The code used to produce the results are available at https://github.com/DENG-MIT/Stiff-PINN.

### 3.1 ROBER problem

We first employed a regular-PINN to solve the ROBER problem. The NN has a single input of time $t$ and outputs the concentrations of three species $[y_1, y_2, y_3]^T$. To avoid the imbalance between the loss of initial conditions and the loss of residuals, we hardcoded the initial conditions to the NN architecture, i.e.,

$$\boldsymbol{y} = \boldsymbol{y_0} + t * NN(\log(t)). \tag{3.1}$$

Therefore, to train the regular-PINN, only the loss functions of the residuals for the three species needed to be minimized. A total of 2500 data points is sampled uniformly in a logarithmic scale of the computational domain (time domain of $t \in [0, 10^5 \, s]$) to evaluate the residuals. The neural network had three hidden layers and 128 nodes per hidden layer, with the activation function of Gaussian Error Linear Units (GELU) [39]. The weights of the NN were initialized using Xavier [40] and optimized via Adam [36], with the default learning rate of 0.001. The training was conducted via mini-batch training with a mini-batch size of 128. The training results are shown in Fig. 4. Overall, the learned regular-PINN can capture species $y_1$ and $y_3$ during the initial stage of $t \in [0, 10 \, s]$, and then substantially deviate from the exact solutions. Although the adaptive weights were implemented for the loss components of the three species following [20,25], the training of regular-PINN for the ROBER problem still failed. The deficiency of adaptive weights in the ROBER problem might be due to the fact that the stiffness here is due to the multiscale nature in the chemical dynamic system, while the stiffness in [25] is attributed to the imbalance among the loss functions of boundary conditions and residuals.

We hypothesize that the failure is due to the physical stiffness of the ROBER problem. A natural way to test this hypothesis is to convert this problem to a non-stiff problem and see whether the regular-PINN would work. We then applied Stiff-PINN to the ROBER problem with QSSA, described by Eq. (2.11), and the results are also shown in Fig. 4. Although the species $y_2$ is not included in the output of the NN, it can be computed via Eq. (2.10). The hyper-parameters for Stiff-PINN are almost identical to those of the regular-PINN except that the Stiff-PINN only outputs two species while the regular-PINN outputs three species. Here, we take the absolute value of $y_1$ and $y_3$ to ensure the validity of Eq. (2.10). As can be seen in Fig. 4, the Stiff-PINN accurately captures all three species profiles, including $y_2$. The history of the loss functions for both regular-PINN and Stiff-PINN are presented in Fig. 5. As expected, the loss function for regular-PINN stays at a very high value, while that of the Stiff-PINN is decreased by a factor of 6 orders of magnitude.

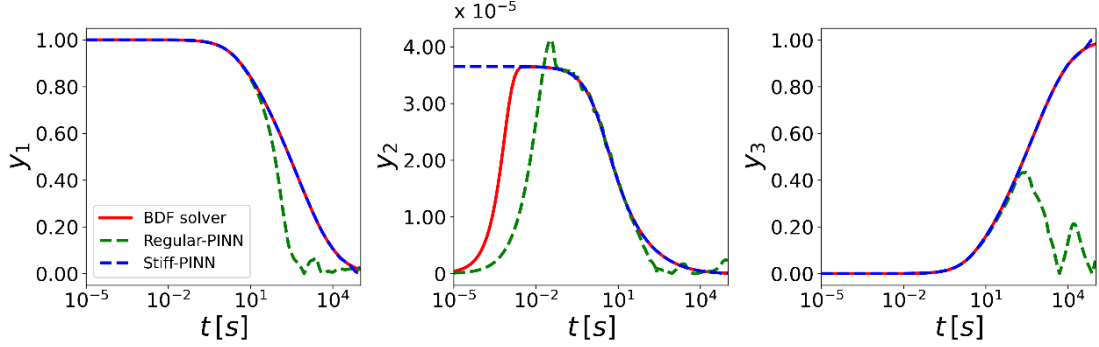

Figure 4. Solutions of the benchmark ROBER problem using the BDF solver (the exact solution), regular-PINN, and Stiff-PINN with QSSA. While the regular-PINN fails to predict the kinetic evolution of the stiff system, Stiff-PINN with QSSA works very well.

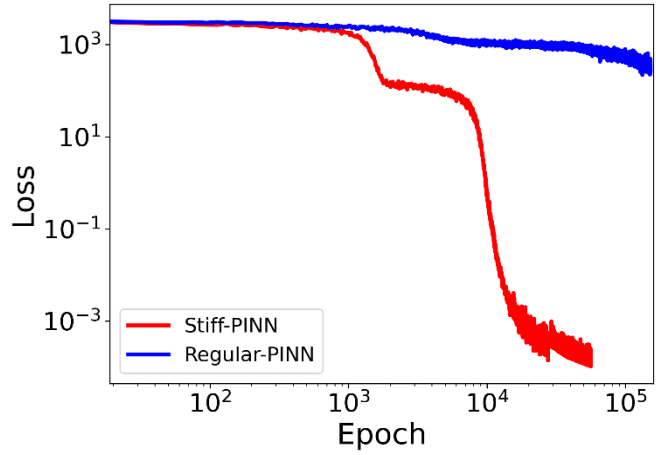

Figure 5. The history of the loss functions of regular-PINN and Stiff-PINN for the ROBER problem. The epoch here corresponds to each parameter update.

While analyzing the gradient flow dynamics is beyond the scope of this work, we draw some intuition on the difference in the performances between Stiff-PINN and regular-PINN. In Stiff-PINN, $y_2$ is eliminated by QSSA, and $y_1, y_3$ are of the same order of magnitude, which makes their evolution much easier to be approximated with a single neural network. In addition, excluding $y_2$ from the loss function also mitigates the imbalance among the residual losses for fast and slow species. However, in regular-PINN, the neural network has to approximate all three species simultaneously, where $y_2 \sim O(10^{-5})$ while $y_1, y_3 \sim O(1)$. The large scale separation is a great challenge to the approximation capacity of the neural network. Therefore, even if the full ODE system is modeled with a sufficiently large and deep neural network, the training would be very challenging.

We shall also briefly discuss the sensitivities of stiff-PINN to the number of nodes and layers of the NN model. The performance is found to be sensitive to the size (neurons × layers) of the NN as shown in Table 1. The general trend is that the performance of deeper neural networks performs better, by comparing the case of 64 × 5 with 64 × 4 and 128 × 3 with 128 × 2, respectively. Empirically, the architecture of three layers and 128 nodes per layer performances best in this work.

Table 1. The root-mean-square-error (RSME) of Stiff-PINN with different NN sizes.

| neurons × layers | $RMSE(y_1)$ | $RMSE(y_3)$ |
| --- | --- | --- |
| 64 × 4 | $4.7954 \times 10^{-3}$ | $7.1774 \times 10^{-3}$ |
| 64 × 5 | $1.4469 \times 10^{-3}$ | $5.0369 \times 10^{-3}$ |
| 128 × 2 | $1.4805 \times 10^{-2}$ | $2.4651 \times 10^{-2}$ |
| 128 × 3 | $1.3628 \times 10^{-3}$ | $2.0699 \times 10^{-3}$ |
| 256 × 1 | $1.5961 \times 10^{-2}$ | $4.5460 \times 10^{-2}$ |

## 3.2 POLLU Problem

The POLLU problem consists of 20 species and 25 reactions. It is an air pollution model developed at The Dutch National Institute of Public Health and Environmental Protection [35]. The POLLU problem can be described mathematically by the following 20 non-linear ODEs shown in Eq. (3.2). Full details about the model can be found in **Error! Reference source not found.**

$$\frac{dy}{dt} = f(y), \ y(0) = y_0, \ y \in R^{20}, \ 0 \leq t \leq 60 \text{ s}. \tag{3.2}$$

As expected, the regular-PINN also failed in the POLLU problem, and the performance is shown in Fig. S1 of the Supporting Information. We then empirically chose 10 species as QSS species based on their maximum concentrations, i.e., lower than 1e-4, and derived the algebraic equations for those 10 species as shown in the Supporting Information. Then, the POLLU problem was reduced to 10 ODEs and 10 algebraic equations from 20 ODEs. The evolution of species concentrations described by the original ODEs and the reduced DAEs by imposing QSSA is shown in Figs. 6 and 7. The QSSA has little effects on the profiles of non-QSS species (Fig. 6), while there are some discrepancies in those of the QSS species (Fig. 7), especially during the initial induction period. However, the assumption does not affect the predictions of the stable species, which are often of greater interest than the unstable intermediates in practice.

Also included in Figs. 6 and 7 are the solutions obtained with the Stiff-PINN approach. To facilitate the training, different fixed weights for each species in the residual loss function were applied to balance them, since the maximum concentrations of the non-QSS species still span several orders of magnitude. Stiff-PINN can accurately predict the evolution of non-QSS species and most of the QSS species. For QSS species $[y_{13}, y_{19}, y_{20}]$, the concentrations of which are close to the non-QSS species, the QSSA tends to induce a larger error compared to the rest of the QSS species, which is expected and consistent with the results using ODE solvers. In general, Stiff-PINN can accurately solve the POLLU problem with stiffness removal via QSSA. We also analyzed the history of the loss functions for regular- and Stiff-PINN, as shown in Fig. 8. Although the loss functions of both regular-PINN and Stiff-PINN decrease as the training goes on, the loss in the Stiff-PINN is four orders of magnitude smaller than that of the regular-PINN.

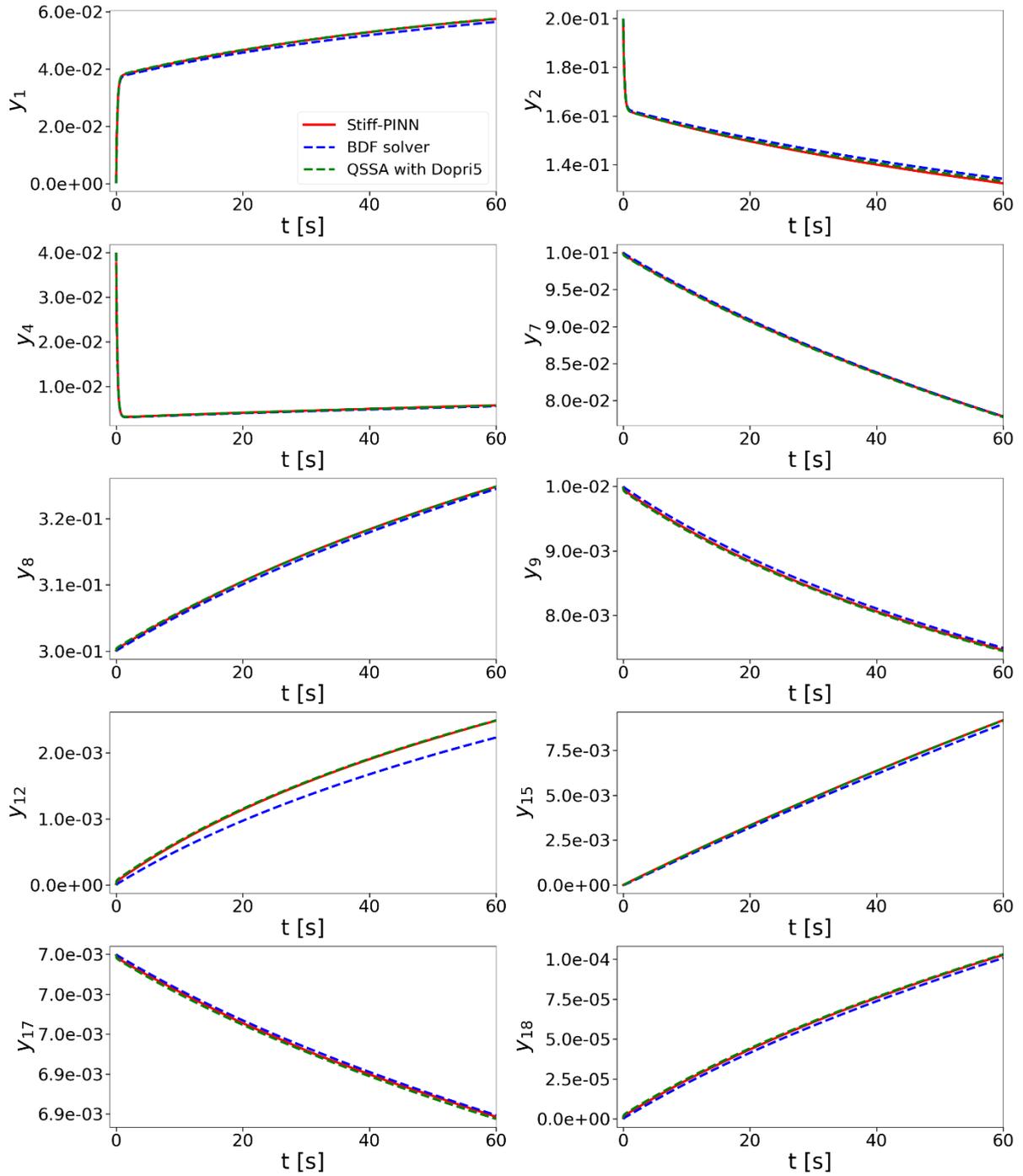

Figure 6. The evolution of the concentration of the ten non-QSS species in the POLLU problem obtained by solving the original 20 ODEs using BDF, the reduced DAEs using Dopri5, and Stiff-PINN.

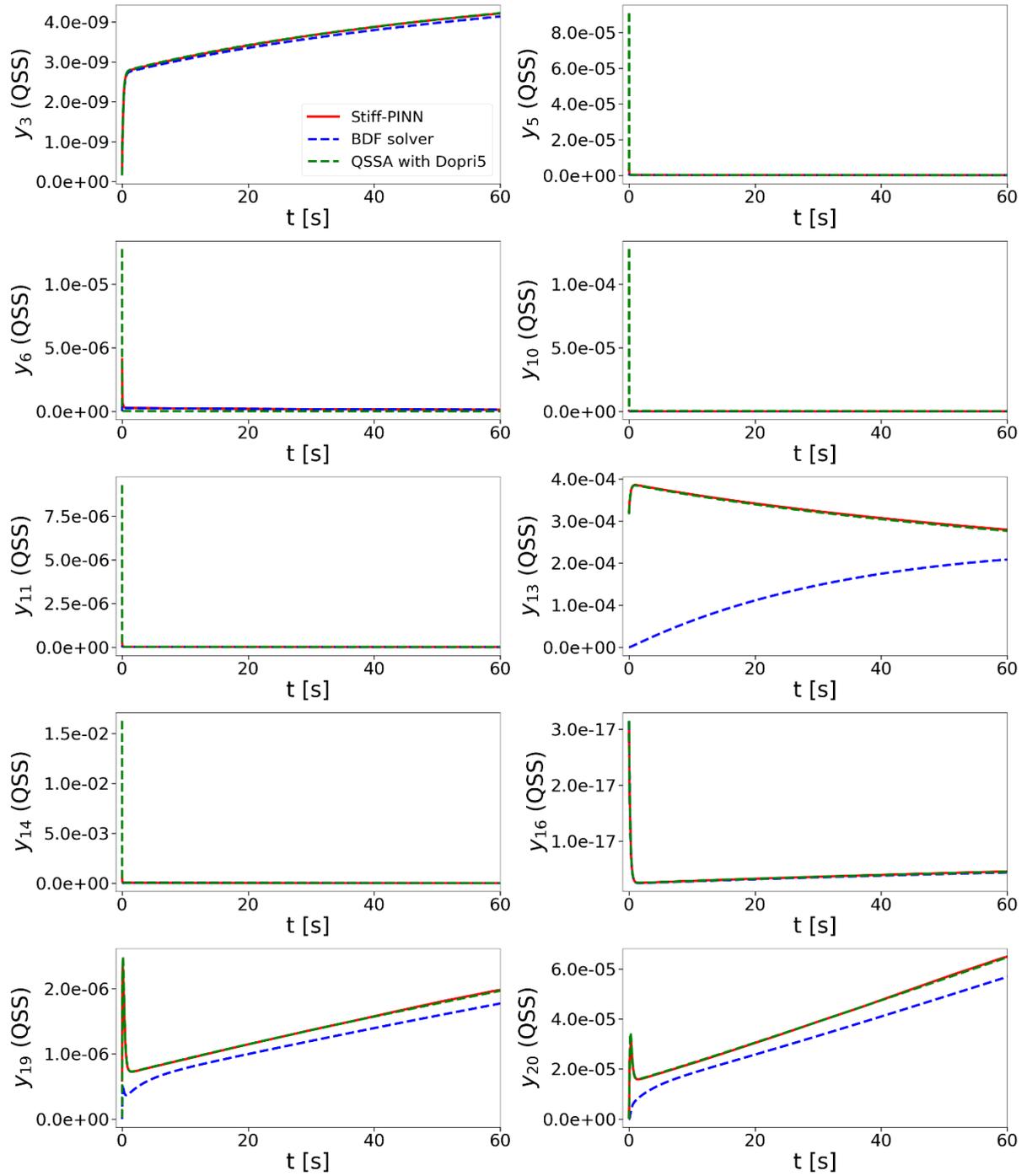

Figure 7. The evolution of the concentration of the ten QSS species in the POLLU problem obtained by solving the original 20 ODEs using BDF, the reduced DAEs using Dopri5, and Stiff-PINN.

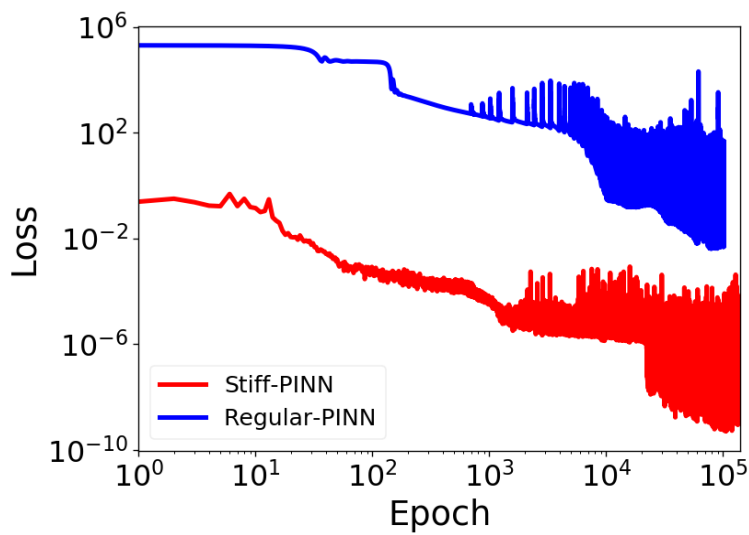

Figure 8. The history of the loss functions of regular-PINN and Stiff-PINN for the POLLU problem.

## 4. Conclusion and Outlook

The performance of the physics-informed neural network in stiff chemical kinetic systems was investigated. In the two classical stiff chemical kinetic systems, ROBER and POLLU, the regular-PINN failed to predict the evolution of the systems. By imposing quasi-steady-state-assumptions on certain species in the kinetic systems and reducing the stiffness, the Stiff-PINN well captured the dynamic responses of the systems. Therefore, it indicates the failure of PINN was due to the physical stiffness, and the proposed Stiff-PINN is an effective and promising approach for stiff chemical kinetic systems.

However, there are still a lot of open questions to be addressed in order to develop a robust and general PINN framework for stiff chemical kinetic problems. We shall highlight some of the challenges learned from performing this work that could guide future advances:

(i) How to handle the fast time scales associated with linear combinations of several species [29,41], rather than single species? More general stiffness removal approaches may help address those challenges, such as computational singular perturbation (CSP) [42,43], intrinsic low dimensional manifolds (ILDM) [44,45], Global Quasi-Linearization (GQL)

$^{29,41}$. Those advanced low-order modeling approaches will also help reduce the approximation error induced by QSSA, e.g., the species of $y_{12}$ in the POLLU problem.

(ii) Manually deriving the QSSA formula for complex chemical systems can be time-consuming and requires a deep understanding of the process. How to identify QSS species and derive the corresponding DAEs automatically to be incorporated into the Stiff-PINN framework? Automatic reduction tools [37,46] and above general low-order approximation approaches (e.g., CSP, ILDM, and GQL) may help tackle the challenge.

(iii) The stiffness in complex chemical systems may not be eliminated, and the reduced system may still show mild stiffness. Advancement in neural network optimizations is required to train PINN for mild stiff systems. A possible solution is exploiting stiff ODE solvers as neural network optimizers [47]. One of the drawbacks of such an approach is that the Hessian matrix of the loss function w.r.t. the neural network parameters are required. Recent advancements in solving stiff ODEs using explicit method [32] and semi-implicit method [33] may help mitigate the requirements of the Hessian matrix. Another possible direction is normalizing the loss functions on-the-fly based on the time scales (e.g., via the eigenvalue of the Jacobian matrix) during the training process.

Finally, while this work is focusing on enabling PINN for stiff systems, the idea of estimating the species that have slow time scales can also help tackle the challenges in other data-driven modeling approaches. For example, it has been shown that training Neural Ordinary Differential Equations for representing kinetic models using neural networks could be challenging for stiff chemical kinetic systems [48], and excluding QSS species from the network could help the training [49].

**Supporting Information.**

Details of the POLLU model including the original full model and the QSSA reduction; Training results for the POLLU model using regular-PINN.

**Acknowledgement**

SD would like to acknowledge the support from the d'Arbeloff Career Development allowance at Massachusetts Institute of Technology.

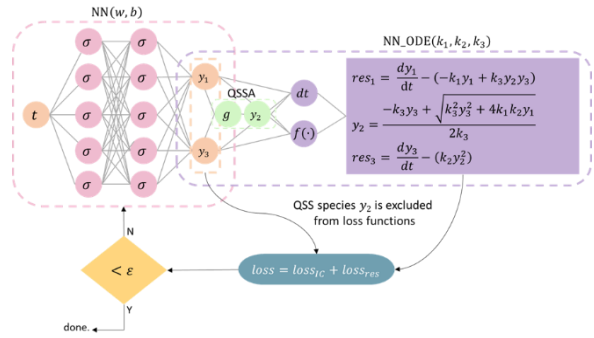

TOC Graphic